# Computational Exploration of Investor Utilities

# Underlying a Portfolio Insurance Strategy

*Dr. M. Khoshnevisan*

*School of Accounting & Finance*

*Griffith University, Australia*

*&*

*Dr. Florentin Smarandache*

*University of New Mexico, USA*

*&*

*Sukanto Bhattacharya*

*School of Information Technology*

*Bond University, Australia*

## Abstract

In this paper we take a look at a simple portfolio insurance strategy using a protective put and computationally derive the investor's governing utility structures underlying such a strategy under alternative market scenarios. Investor utility is deemed to increase with an increase in the excess equity generated by the portfolio insurance strategy over a simple investment strategy without any insurance. Three alternative market scenarios (probability spaces) have been explored – "Down", "Neutral" and "Up", categorized according to whether the price of the underlying security is most likely to go down, stay unchanged or go up. The methodology used is computational, primarily based on simulation and numerical extrapolation. The Arrow-Pratt measure of risk aversion has been used to determine how the investors react towards risk under the different scenarios.

**Key words**

Option pricing, investment risk, portfolio insurance, utility theory, behavioral economics

**Introduction:**

Basically, a derivative financial asset is a legal contract between two parties – a buyer and a seller, whereby the former receives a rightful claim on an underlying asset while the latter has the corresponding liability of making good that claim, in exchange for a mutually agreed consideration. While many derivative securities are traded on the floors of exchanges just like ordinary securities, some derivatives are not exchange-traded at all. These are called OTC (Over-the-Counter) derivatives, which are contracts not traded on organized exchanges but rather negotiated privately between parties and are especially tailor-made to suit the nature of the underlying assets and the pay-offs desired therefrom. While countless papers have been written on the mathematics of option pricing formulation, surprisingly little work has been done in the area of exploring the exact nature of investor utility structures that underlie investment in derivative financial assets. This is an area we deem to be of tremendous interest both from the point of view of mainstream financial economics as well as from the point of view of a more recent and more esoteric perspective of *behavioral economics*.

**The basic building blocks of derivative assets:**

Forward Contract

A contract to buy or sell a specified amount of a designated commodity, currency, security, or financial instrument at a known date in the future and at a price set at the time the contract is made. Forward contracts are negotiated between the contracting parties and are not traded on organized exchanges.

Futures Contract

Quite similar to a forwards contract – this is a contract to buy or sell a specified amount of a designated commodity, currency, security, or financial instrument at a known date in the future and at a price set at the time the contract is made. What primarily distinguishes forward contracts from futures contracts is that the latter are traded on organized exchanges and are thus standardized. These contracts are marked to market daily, with profits and losses settled in cash at the end of the trading day.

Swap Contract

A private contract between two parties to exchange cash flows in the future according to some prearranged formula. The most common type of swap is the "plain vanilla" interest rate swap, in which the first party agrees to pay the second party cash flows equal to interest at a predetermined f*ixed* rate on a notional principal. The second party agrees to pay the first party cash flows equal to interest at a *floating* rate on the same notional principal. Both payment streams are denominated in the same currency. Another common type of swap is the currency swap. This contract calls for the counter-parties to exchange specific amounts of two different currencies at the outset, which are repaid over time according to a prearranged formula that reflects amortization and interest payments.

Option Contract

A contract that gives its owner the right, but not the obligation, to buy or sell a specified asset at a stipulated price, called the strike price. Contracts that give owners the right to buy are referred to as *call* options and contracts that give the owner the right to sell are called *put* options. Options include both standardized products that trade on organized exchanges and customized contracts between private parties.

In our present analysis we will be restricted exclusively to portfolio insurance strategy using a long position in put options and explore the utility structures derivable therefrom.

The simplest option contracts (also called *plain vanilla options*) are of two basic types – **call** and **put**. The call option is a right to buy (or call up) some underlying asset at or within a specific future date for a specific price called the strike price. The put option is a right to sell (or put through) some underlying asset at or within a specified date – again for a pre-determined strike price. The options come with no obligations attached – it is totally the discretion of the option holder to decide whether or not to exercise the same.

The pay-off function (from an option buyer's viewpoint) emanating from a call option is given as $P_{call} = Max\,[(S_T - X),\,0]$. Here, $S_T$ is the price of the underlying asset on maturity and X is the strike price of the option. Similarly, for a put option, the pay-off function is given as $P_{put} = Max\ [(X – S_T),\ 0]$. The implicit assumption in this case is that the options can only be exercised on the maturity date and not earlier. Such options are called **European options**. If the holder of an option contract is allowed to exercise the same any time on or before the day of maturity, it is termed an **American option**. A third, not-so-common category is one where the holder can exercise the option only on specified dates prior to its maturity. These are termed **Bermudan options**. The options we refer to in this paper will all be European type only but methodological extensions are possible to extend our analysis to also include American or even Bermudan options.

**Investor's utility structures governing the purchase of plain vanilla option contracts:**

Let us assume that an underlying asset priced at S at time t will go up or down by Δs or stay unchanged at time T either with probabilities $p_U(u)$, $p_U(d)$ and $p_U(n)$ respectively contingent upon the occurrence of event U, or with probabilities $p_D(u)$, $p_D(d)$ and $p_D(n)$ respectively contingent upon the occurrence of event D, or with probabilities $p_N(u)$, $p_N(d)$ and $p_N(n)$ respectively contingent upon the occurrence of event N, in the time period (T – t). This, by the way, is comparable to the analytical framework that is exploited in option pricing using the numerical method of **trinomial trees**. The trinomial tree algorithm is mainly used in the pricing of the non-European options where no closed-form pricing formula exists.

**Theorem:**

Let $P_U$, $P_D$ and $P_N$ be the three probability distributions contingent upon events U, D and N respectively. Then we have *a consistent preference relation* for a call buyer such that ***$P_U$ is strictly preferred to $P_N$ and $P_N$ is strictly preferred to $P_D$*** and *a corresponding consistent preference relation* for a put buyer such that ***$P_D$ is strictly preferred to $P_N$ and $P_N$ is strictly preferred to $P_U$***.

**Proof:**

Case I: Investor buys a call option for \$C maturing at time T having a strike price of \$X on the underlying asset. We modify the call pay-off function slightly such that we now have the pay-off function as: ***$P_{call} = Max\ (S_T - X - C_{price},\ - C_{price})$***.

Event U: $E_U\ (Call) = [(S + e^{-r(T-t)} \Delta s)\ p_U\ (u) + (S - e^{-r(T-t)} \Delta s)\ p_U\ (d) + S\ p_U\ (n)] - C - Xe^{-r(T-t)}$

$$= [S + e^{-r(T-t)} \Delta s\ \{p_U\ (u) - p_U\ (d)\}] - C - Xe^{-r(T-t)} \ldots p_U\ (u) > p_U\ (d)$$

Therefore, $E\ (P_{call}) = Max\ [S + e^{-r(T-t)} \{\Delta s\ (p_U\ (u) - p_U\ (d)) - X\} - C, - C]$ **… (1)**

Event D: $E_D\ (Call) = [(S + e^{-r(T-t)} \Delta s)\ p_D\ (u) + (S - e^{-r(T-t)} \Delta s)\ p_D\ (d) + S\ p_D\ (n)] - C - Xe^{-r(T-t)}$

$$= [S + e^{-r(T-t)} \Delta s\ \{p_D\ (u) - p_D\ (d)\}] - C - Xe^{-r(T-t)} \ldots p_D\ (u) < p_D\ (d)$$

Therefore, $E\ (P_{call}) = Max\ [S - e^{-r(T-t)} \{\Delta s\ (p_D\ (d) - p_D\ (u)) + X\} - C, - C]$ **… (2)**

Event N: $E_N$ (Call) = $[(S + e^{-r(T-t)} \Delta s) p_N(u) + (S - e^{-r(T-t)} \Delta s) p_N(d) + S p_N(n)] - C - Xe^{-r(T-t)}$

$$= [S + e^{-r(T-t)} \Delta s \{p_N(u) - p_N(d)\}] - C - Xe^{-r(T-t)}$$

$$= S - C - Xe^{-r(T-t)} \ldots p_N(u) = p_N(d)$$

Therefore, $E(P_{call}) = \text{Max}[S - Xe^{-r(T-t)} - C, -C]$ ... (3)

Case II: Investor buys a put option for $P maturing at time T having a strike price of $X on the underlying asset. Again we modify the pay-off function such that we now have the pay-off function as: $\boldsymbol{P_{put} = Max(X - S_T - P_{price}, -P_{price})}$.

Event U: $E_U$ (Put) = $Xe^{-r(T-t)} - [\{(S + e^{-r(T-t)} \Delta s) p_U(u) + (S - e^{-r(T-t)} \Delta s) p_U(d) + S p_U(n)\} + P]$

$$= Xe^{-r(T-t)} - [S + e^{-r(T-t)} \Delta s \{p_U(u) - p_U(d)\} + P]$$

$$= Xe^{-r(T-t)} - [S + e^{-r(T-t)} \Delta s \{p_U(u) - p_U(d)\} + (C + Xe^{-r(T-t)} - S)] \ldots \textit{put-call parity}$$

$$= -e^{-r(T-t)} \Delta s \{p_U(u) - p_U(d)\} - C$$

Therefore, $E(P_{put}) = \text{Max}[-e^{-r(T-t)} \Delta s \{p_U(u) - p_U(d)\} - C, -P]$

$$= \text{Max}[-e^{-r(T-t)} \Delta s \{p_U(u) - p_U(d)\} - C, -(C + Xe^{-r(T-t)} - S)] \ldots \textbf{(4)}$$

Event D: $E_D$ (Put) = $Xe^{-r(T-t)} - [\{(S + e^{-r(T-t)} \Delta s) p_D(u) + (S - e^{-r(T-t)} \Delta s) p_D(d) + S p_D(n)\} + P]$

$$= Xe^{-r(T-t)} - [S + e^{-r(T-t)} \Delta s \{p_D(u) - p_D(d)\} + P]$$

$$= Xe^{-r(T-t)} - [S + e^{-r(T-t)} \Delta s \{p_U(u) - p_U(d)\} + (C + Xe^{-r(T-t)} - S)] \ldots \textit{put-call parity}$$

$$= e^{-r(T-t)} \Delta s \{p_D(d) - p_D(u)\} - C$$

Therefore, E ($P_{put}$) = Max [$e^{-r(T-t)}$ Δs {$p_D$ (d) – $p_D$ (u)} – C, – P]

$$= \text{Max} [e^{-r(T-t)} \Delta s \{p_D (d) - p_D (u)\} - C, - (C + Xe^{-r(T-t)} - S)] \quad \ldots (5)$$

Event N: $E_N$ (Put) = $Xe^{-r(T-t)}$ – [{(S + $e^{-r(T-t)}$ Δs) $p_N$ (u) + (S – $e^{-r(T-t)}$ Δs) $p_N$ (d) + S $p_N$ (n)} + P]

$$= Xe^{-r(T-t)} - [S + e^{-r(T-t)} \Delta s \{p_N (u) - p_N (d)\} + P]$$

$$= Xe^{-r(T-t)} - (S + P)$$

$$= (Xe^{-r(T-t)} - S) - \{C + (Xe^{-r(T-t)} - S)\} \ldots \textit{put-call parity}$$

$$= - C$$

Therefore, E ($P_{put}$) = Max [– C, – P]

$$= \text{Max} [-C, - (C + Xe^{-r(T-t)} - S)] \quad \ldots (6)$$

From equations (4), (5) and (6) we see that **$E_U$ (Put) < $E_N$ (Put) < $E_D$ (Put)** and hence it is proved why we have the *consistent preference relation $P_D$ is strictly preferred to $P_N$ and $P_N$ is strictly preferred to $P_U$* from a put buyer's point of view. The call buyer's consistent preference relation is also explainable likewise.

We can now proceed to computationally derive the associated utility structures using a Monte Carlo discrete-event simulation approach to estimate the change in equity following a particular investment strategy under each of the aforementioned event spaces.

**Computational derivation of investor's utility curves under a *protective put* strategy:**

There is a more or less well-established theory of utility maximization in case of deductible insurance policy on non-financial assets whereby the basic underlying assumption is that cost of insurance is a convex function of the expected indemnification. Such an assumption has been showed to satisfy the sufficiency condition for expected utility maximization when individual preferences exhibit risk aversion. The final wealth function at end of the insurance period is given as follows:

$$Z_T = Z_0 + M - x + I(x) - C(D) \qquad \ldots (7)$$

Here $Z_T$ is the final wealth at time $t = T$, $Z_0$ is the initial wealth at time $t = 0$, x is a random loss variable, I (x) is the indemnification function, C (x) is the cost of insurance and $0 \leq D \leq M$ is the level of the deductible. However the parallels that can be drawn between ordinary insurance and portfolio insurance is different when the portfolio consists of financial assets being continuously traded on the floors of organized financial markets. While the form of an insurance contract might look familiar – an assured value in return for a price – the mechanism of providing such assurance will have to be quite different because unlike other tangible assets like houses or cars, when one portfolio of financial assets gets knocked down, virtually all others are likely to follow suit making "risk pooling", the typical method of insurance, quite inadequate for portfolio insurance. Derivative assets like options do provide a suitable mechanism for portfolio insurance.

If the market is likely to move adversely, holding a long put alongside ensures that the investor is better off than just holding a long position in the underlying asset. The long put offers the investor some kind of *price insurance* in case the market goes down. This strategy is known in derivatives parlance as a *protective put*. The strategy effectively puts a floor on the downside deviations without cutting off the upside by too much. From the expected changes in investor's equity we can computationally derive his or her utility curves under the strategies $A_1$ and $A_2$ in each of the three probability spaces D, N and U.

The following hypothetical data have been assumed to calculate the simulated put price:

S = $50.00 (purchase price of the underlying security)

X = $55.00 (put strike price)

(T – t) = 1 (single period investment horizon)

Risk-free rate = 5%

The put options has been valued by Monte Carlo simulation of a trinomial tree using a customized MS-Excel spreadsheet for one hundred independent replications in each case.

Event space: D Strategy: $A_1$ (Long underlying asset)

Instance (i): $(-)\Delta S = \$5.00$, $(+)\Delta S = \$15.00$

Table 1

| Price movement | Probability | Expected Δ Equity |
|---|---|---|
| Up (+ $15.00) | 0.1 | $1.50 |
| Neutral ($0.00) | 0.3 | $0.00 |
| Down (– $5.00) | 0.6 | ($3.00) |
|  |  | Σ = ($1.50) |

To see how the expected change in investor's equity goes up with an increased upside potential we will double the possible up movement at each of the next two stages while keeping the down movement unaltered. This should enable us to account for any possible loss of investor utility by way of the cost of using a portfolio insurance strategy.

Instance (ii): $(+)\Delta S = \$30.00$

Table 2

| Price movement | Probability | Expected Δ Equity |
|---|---|---|
| Up (+ $30.00) | 0.1 | $3.00 |
| Neutral ($0.00) | 0.3 | $0.00 |
| Down (− $5.00) | 0.6 | ($3.00) |
| | | $\Sigma = \$0.00$ |

Instance (iii): $(+)\Delta S = \$60.00$

Table 3

| Price movement | Probability | Expected Δ Equity |
|---|---|---|
| Up (+ $60.00) | 0.1 | $6.00 |
| Neutral ($0.00) | 0.3 | $0.00 |
| Down (− $5.00) | 0.6 | ($3.00) |
| | | $\Sigma = \$3.00$ |

Event space: D Strategy: $A_2$ (Long underlying asset + long put)

Instance (i): $(-)\Delta S = \$5.00$, $(+)\Delta S = \$15.00$

Table 4

| Simulated put price | $6.99 |
|---|---|
| Variance | $11.63 |
| Simulated asset value | $48.95 |
| Variance | $43.58 |

**Table 5**

| Price movement | Probability | Expected Δ Equity | Expected excess equity | Utility index |
|---|---|---|---|---|
| Up (+ $8.01) | 0.1 | $0.801 | | |
| Neutral (– $1.99) | 0.3 | ($0.597) | | |
| Down (– $1.99) | 0.6 | ($1.194) | | |
| | | Σ = (–$0.99) | $0.51 | ≈ 0.333 |

Instance (ii): (+)ΔS = $30.00

**Table 6**

| | |
|---|---|
| Simulated put price | $6.75 |
| Variance | $13.33 |
| Simulated asset value | $52.15 |
| Variance | $164.78 |

**Table 7**

| Price movement | Probability | Expected Δ Equity | Expected excess equity | Utility index |
|---|---|---|---|---|
| Up (+ $23.25) | 0.1 | $2.325 | | |
| Neutral (– $1.75) | 0.3 | ($0.525) | | |
| Down (– $1.75) | 0.6 | ($1.05) | | |
| | | Σ = $0.75 | $0.75 | ≈ 0.666 |

Instance (iii): (+)ΔS = $60.00

**Table 8**

| | |
|---|---|
| Simulated put price | $6.71 |
| Variance | $12.38 |
| Simulated asset value | $56.20 |
| Variance | $520.77 |

Table 10

| Price movement | Probability | Expected Δ Equity | Expected excess equity | Utility index |
|---|---|---|---|---|
| Up (+ $53.29) | 0.1 | $5.329 | | |
| Neutral (– $1.71) | 0.3 | ($0.513) | | |
| Down (– $1.71) | 0.6 | ($1.026) | | |
| | | Σ = $3.79 | $0.79 | ≈ 0.999 |

**Figure 1**

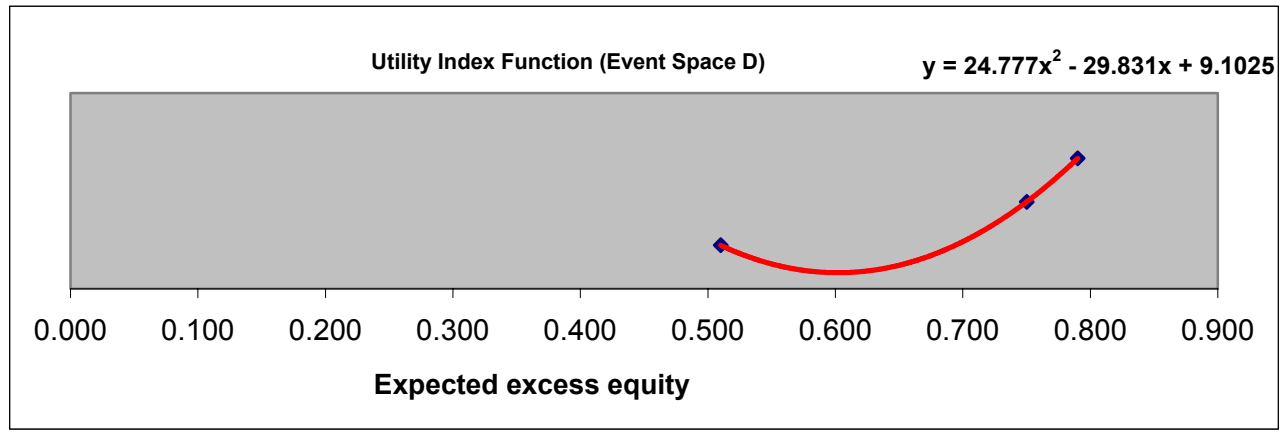

The utility function as obtained above is *convex* in probability space D, which indicates that the protective strategy can make the investor risk-loving even when the market is expected to move in an adverse direction, as the expected payoff from the put option largely neutralizes the likely erosion of security value at an affordable insurance cost! This seems in line with intuitive behavioral reasoning, as investors with a viable downside protection will become more aggressive in their approach than they would be without it implying markedly lowered risk averseness for the investors with insurance.

Event space: N Strategy: $A_1$ (Long underlying asset)

Instance (i): $(-)\Delta S = \$5.00$, $(+)\Delta S = \$15.00$

**Table 11**

| Price movement | Probability | Expected Δ Equity |
|---|---|---|
| Up (+ $15.00) | 0.2 | $3.00 |
| Neutral ($0.00) | 0.6 | $0.00 |
| Down (− $5.00) | 0.2 | ($1.00) <br> Σ = $2.00 |

Instance (ii):  (+)ΔS = $30.00

**Table 12**

| Price movement | Probability | Expected Δ Equity |
|---|---|---|
| Up (+ $30.00) | 0.2 | $6.00 |
| Neutral ($0.00) | 0.6 | $0.00 |
| Down (− $5.00) | 0.2 | ($1.00) <br> Σ = $5.00 |

Instance (iii):  (+)ΔS = $60.00

**Table 13**

| Price movement | Probability | Expected Δ Equity |
|---|---|---|
| Up (+ $60.00) | 0.2 | $12.00 |
| Neutral ($0.00) | 0.6 | $0.00 |
| Down (− $5) | 0.2 | ($1.00) <br> Σ = $11.00 |

Event space: N Strategy: $A_2$ (Long underlying asset + long put)

Instance (i): (−)ΔS = $5.00, (+)ΔS = $15.00

**Table 14**

| | |
|---|---|
| Simulated put price | $4.85 |
| Variance | $9.59 |
| Simulated asset value | $51.90 |
| Variance | $47.36 |

**Table 15**

| Price movement | Probability | Expected Δ Equity | Expected excess equity | Utility index |
|---|---|---|---|---|
| Up (+ $11.15) | 0.2 | $2.23 | | |
| Neutral (+ $0.15) | 0.6 | $0.09 | | |
| Down (+ $0.15) | 0.2 | $0.03 | | |
| | | Σ = $2.35 | $0.35 | ≈ 0.999 |

Instance (ii): (+)ΔS = $30.00

**Table 16**

| | |
|---|---|
| Simulated put price | $4.80 |
| Variance | $9.82 |
| Simulated asset value | $55.20 |
| Variance | $169.15 |

**Table 17**

| Price movement | Probability | Expected Δ Equity | Expected excess equity | Utility index |
|---|---|---|---|---|
| Up (+ $25.20) | 0.2 | $5.04 | | |
| Neutral (+ $0.20) | 0.6 | $0.12 | | |
| Down (+ $0.20) | 0.2 | $0.04 | | |
| | | Σ = $5.20 | $0.20 | ≈ 0.333 |

Instance (iii): (+)ΔS = $60.00

**Table 18**

| | |
|---|---|
| Simulated put price | $4.76 |
| Variance | $8.68 |
| Simulated asset value | $60.45 |
| Variance | $585.40 |

**Table 19**

| Price movement | Probability | Expected Δ Equity | Expected excess equity | Utility index |
|---|---|---|---|---|
| Up (+ $55.24) | 0.2 | $11.048 | | |
| Neutral (+ $0.24) | 0.6 | $0.144 | | |
| Down (+ $0.24) | 0.2 | $0.048 | | |
| | | Σ = $11.24 | $0.24 | ≈ 0.666 |

**Figure 2**

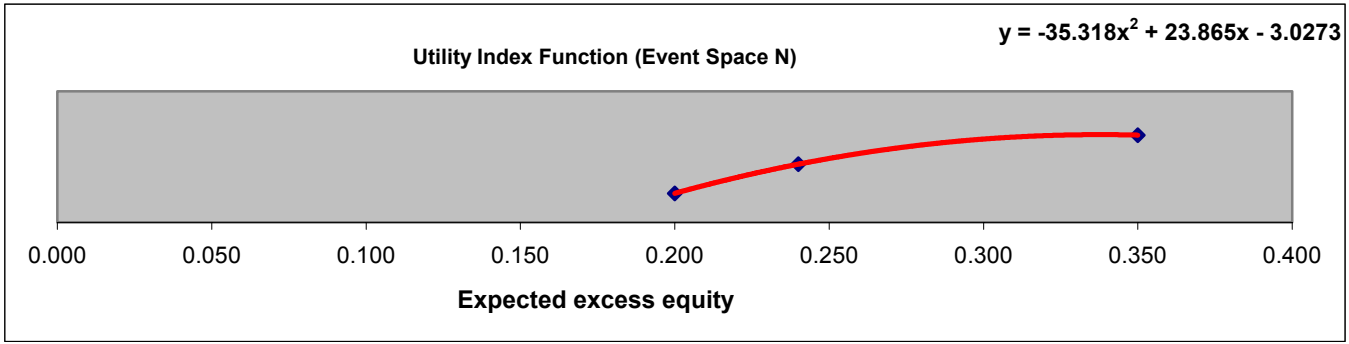

The utility function as obtained above is *concave* in probability space N, which indicates that the insurance provided by the protective strategy can no longer make the investor risk-loving as the expected value of the insurance is offset by the cost of buying the put! This is again in line with intuitive behavioral reasoning because if the market is equally likely to move up or down and more likely to stay unmoved the investor would deem himself or herself better off not buying the insurance because in order to have the insurance i.e. the put option it is necessary to pay an out-of-pocket cost, which may not be offset by the expected payoff from the put option under the prevalent market scenario.

Event space: U Strategy: $A_1$ (Long underlying asset)

Instance (i):  $(-)\Delta S = \$5.00$, $(+)\Delta S = \$15.00$

**Table 20**

| Price movement | Probability | Expected Δ Equity |
|---|---|---|
| Up (+ $15.00) | 0.6 | $9.00 |
| Neutral ($0.00) | 0.3 | $0.00 |
| Down (– $5.00) | 0.1 | ($0.50) |
|  |  | Σ = $8.50 |

Instance (ii): $(+)\Delta S = \$30.00$

**Table 21**

| Price movement | Probability | Expected Δ Equity |
|---|---|---|
| Up (+ $30.00) | 0.6 | $18.00 |
| Neutral ($0.00) | 0.3 | $0.00 |
| Down (− $5.00) | 0.1 | ($0.50) |
| | | Σ = $17.50 |

Instance (iii): $(+)\Delta S = \$60.00$

**Table 22**

| Price movement | Probability | Expected Δ Equity |
|---|---|---|
| Up (+ $60.00) | 0.6 | $36.00 |
| Neutral ($0.00) | 0.3 | $0.00 |
| Down (− $5) | 0.1 | ($0.50) |
| | | Σ = $35.50 |

Event space: U Strategy: $A_2$ (Long underlying asset + long put)

Instance (i): $(-)\Delta S = \$5.00$, $(+)\Delta S = \$15.00$

**Table 23**

| Simulated put price | $2.28 |
|---|---|
| Variance | $9.36 |
| Simulated asset value | $58.60 |
| Variance | $63.68 |

**Table 24**

| Price movement | Probability | Expected Δ Equity | Expected excess equity | Utility index |
|---|---|---|---|---|
| Up (+ $12.72) | 0.6 | $7.632 | | |
| Neutral (+ $2.72) | 0.3 | $0.816 | | |
| Down (+ $2.72) | 0.1 | $0.272 | | |
| | | Σ = $8.72 | $0.22 | ≈ 0.333 |

Instance (ii): $(+)\Delta S = \$30.00$

**Table 25**

| Simulated put price | $2.14 |
|---|---|
| Variance | $10.23 |
| Simulated asset value | $69.00 |
| Variance | $228.79 |

**Table 26**

| Price movement | Probability | Expected Δ Equity | Expected excess equity | Utility index |
|---|---|---|---|---|
| Up (+ $27.86) | 0.6 | $16.716 | | |
| Neutral (+ $2.86) | 0.3 | $0.858 | | |
| Down (+ $2.86) | 0.1 | $0.286 | | |
| | | Σ = $17.86 | $0.36 | ≈ 0.666 |

Instance (iii): $(+)\Delta S = \$60.00$

**Table 27**

| Simulated put price | $2.09 |
|---|---|
| Variance | $9.74 |
| Simulated asset value | $88.55 |
| Variance | $864.80 |

**Table 28**

| Price movement | Probability | Expected Δ Equity | Expected excess equity | Utility index |
|---|---|---|---|---|
| Up (+ $57.91) | 0.6 | $34.746 | | |
| Neutral (+ $2.91) | 0.3 | $0.873 | | |
| Down (+ $2.91) | 0.1 | $0.291 | | |
| | | Σ = $35.91 | $0.41 | ≈ 0.999 |

**Figure 3**

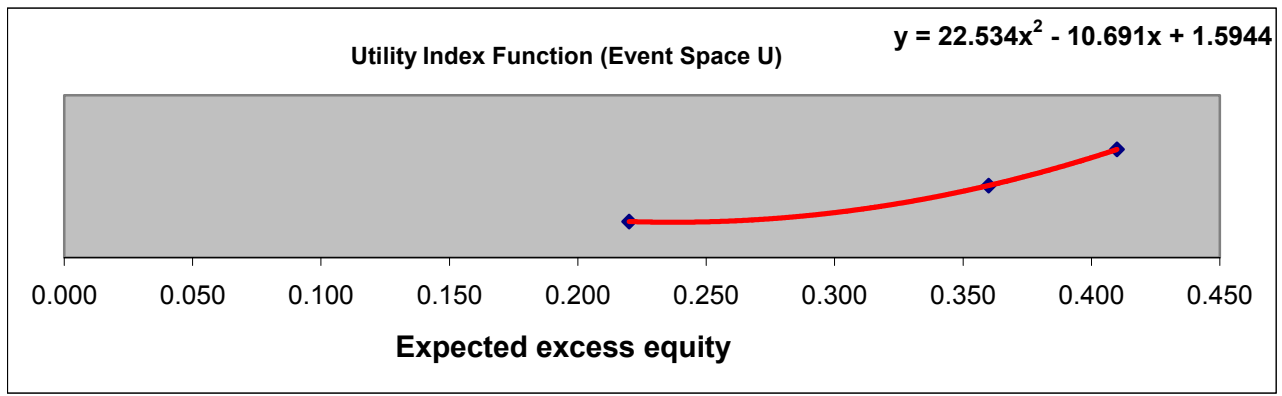

In accordance with intuitive, behavioral reasoning the utility function is again seen to be *convex* in the probability space U, which is probably attributable to the fact that while the market is expected to move in a favourable direction the put option nevertheless keeps the downside protected while costing less than the expected payoff on exercise thereby fostering a risk-loving attitude in the investor as he gets to enjoy the best of both worlds.

**Note:** Particular values assigned to the utility indices won't affect the essential mathematical structure of the utility curve – but only cause a scale shift in the parameters. For example, the indices could easily have been taken as (0.111, 0.555, 0.999) - these assigned values should not have any computational significance as long as all they all lie within the conventional interval (0, 1]. Repeated simulations have shown that the investor would be considered extremely unlucky to get an excess return less than the minimum excess return obtained or extremely lucky to get an excess return more than the maximum excess return obtained under each of the event spaces. Hence, the maximum and minimum expected excess equity within a particular event space should correspond to the lowest and highest utility indices and the utility derived from the median excess equity should then naturally occupy the middle position. As long as this is the case, there will be

no alteration in the fundamental mathematical structure of the investor's utility functions no matter what index values are assigned to his or her utility from expected excess equity.

**Extrapolating the ranges of investor's risk aversion within each probability space:**

For a continuous, twice-differentiable utility function u (x), the *Arrow-Pratt measure of absolute risk aversion* **(ARA)** is given as follows:

$$\lambda (x) = -[d^2u(x)/dx^2][du(x)/dx]^{-1} \qquad \ldots (8)$$

$\lambda (x) > 0$ if u is monotonically increasing and *strictly concave* as in case of a risk-averse investor having u'' (x) < 0. Obviously, $\lambda (x) = 0$ for the risk-neutral investor with a *linear* utility function having u'' (x) = 0 while $\lambda (x) < 0$ for the risk-loving investor with a *strictly convex* utility function having u'' (x) > 0.

Case I: Probability Space D:

u (x) = $24.777x^2$ − 29.831x + 9.1025, u' (x) = 49.554x − 29.831 and u''(x) = 49.554. Thus $\lambda (x)$ = − 49.554/(49.554x − 29.831). Therefore, given the *convex utility function*, the defining range is $\lambda (x) < 0$ i.e. (49.554x − 29.831) < 0 or x < 0.60199.

Case II: Probability Space N:

u (x) = -35.318 $x^2$ + 23.865x − 3.0273, u' (x) = -70.636x + 23.865 and u''(x) = −70.636.

Thus, $\lambda(x) = -[-70.636/(-70.636x + 23.865)] = 70.636/(-70.636x + 23.865)$. Therefore, given the *concave utility function*, the defining range is $\lambda(x) > 0$, i.e. we have the denominator $(-70.636x + 23.865) > 0$ or $x > 0.33786$.

Case III: Probability Space U:

$u(x) = 22.534x^2 - 10.691x + 1.5944$, $u'(x) = 45.068x - 10.691$ and $u''(x) = 45.068$. Thus $\lambda(x) = -45.068/(45.068x - 10.691)$. Therefore, given the *convex utility function*, the defining range is $\lambda(x) < 0$ i.e. $(45.068x - 10.691) < 0$ or $x < 0.23722$.

These defining ranges as evaluated above will however depend on the parameters of the utility function and will therefore be different for different investors according to the values assigned to his or her utility indices corresponding to the expected excess equity.

In general, if we have a parabolic utility function $u(x) = a + bx - cx^2$, where $c > 0$ ensures concavity, then we have $u'(x) = b - 2cx$ and $u''(x) = -2c$. The Arrow-Pratt measure is given by $\lambda(x) = 2c/(b-2cx)$. Therefore, for $\lambda(x) \geq 0$, we need $b \geq 2cx$, thus it can only apply for a limited range of x. Notice that $\lambda'(x) \geq 0$ up to where $x = b/2c$. Beyond that, marginal utility is negative - i.e. beyond this level of equity, utility *declines*. One more implication is that there is an increasing apparent unwillingness to take risk as their equity increases, i.e. with larger excess equity investors are less willing to take risks as concave, parabolic utility functions exhibit increasing absolute risk aversion (IARA).

People sometimes use a past outcome as a critical factor in evaluating the likely outcome from a risky decision to be taken in the present. Also it has been experimentally demonstrated that decisions can be taken in violation of conditional preference relations. This has been the crux of a whole body of behavioral utility theory developed on the basis of what has come to be known as *non-expected utility* following the famous work in *prospect theory* (Kahneman and Tversky, 1979). It has been empirically demonstrated that people are willing to take more risks immediately following gains and take less risks immediately following losses with the probability distribution of the payoffs remaining unchanged. Also decisions are affected more by instantaneous utility resulting from immediate gains than by disutility resulting from the cumulative magnitude of likely losses as in the assessment of health risks from addictive alcohol consumption. It has also been seen in experimental psychology studies that generated explanations cause a greater degree of belief persistence than provided explanations. This is due to a psychological miscalibration whereby people tend to be guided by outcomes in their most recent memory. In the face of all these challenges to the expected utility paradigm, it must however be noted that the utility structures underlying the behavior of investors with loss insurance in the three different market scenarios as derived above are *independent of any psychological miscalibration* on the part of the individual based on prior history of positive or negative payoffs but rather are a direct statistical consequence of the portfolio insurance strategy itself and the expected payoffs likely to follow from such a strategy.

**Conclusion:**

In this paper we have computationally examined the implications on investor's utility of a simple option strategy of portfolio insurance under alternative market scenarios and is hence novel both in content as well as context. We have found that such insurance strategies can indeed have quite interesting governing utility structures underlying them. The expected excess payoffs from an insurance strategy can make the investor risk-loving when it is known with a relatively high prior probability that the market will either move in an adverse direction or in a favourable direction. The investor seems to display risk-averseness only when the market is equally likely to move in either direction and has a relatively high prior probability of staying unmoved. The door is now open for further research along these lines going deep into the governing utility structures that may underlie more complex derivative trading strategies and portfolio insurance schemes.

*******